\documentclass[letterpaper, 10 pt, conference]{ieeeconf}  % Comment this line out if you need a4paper

\IEEEoverridecommandlockouts                              % This command is only needed if
\usepackage{lmodern} % Enhanced version of Computer Modern font family
\usepackage[utf8]{inputenc} % Allows input of UTF-8 characters
\usepackage{mathtools} % Extends amsmath with additional features
\usepackage{amsfonts} % Additional math fonts like \mathbb
\usepackage{amssymb} % Additional math symbols like \mathbb and \mathfrak
\usepackage{physics} % Commands for typesetting physics
\usepackage{graphicx} % Include graphics
\usepackage{subcaption} % Sub-captions for sub-figures and sub-tables
\usepackage{algorithm} % Floating environment for algorithms
\usepackage[noend]{algpseudocode} % Pseudocode environment for algorithms
\usepackage[hidelinks]{hyperref} % Adds hyperlinks without visual effects
\usepackage{xspace} % Ensures correct spacing after commands
\usepackage{bm} % Acces bold symbols in math mode
\usepackage{dsfont} % Provides double struck math fonts
\usepackage{xcolor} % Colors
\usepackage{tikz} % TikZ for creating graphics
\usepackage{pgfplots} % PGFPlots for plotting
\pgfplotsset{compat=1.18}

\newcommand{\ie}{\textit{i.e.,}\xspace} %
\newcommand{\eg}{\textit{e.g.,}\xspace} %

\newcommand\Acal{{\mathcal{A}}}
\newcommand\Pcal{{\mathcal{P}}}
\newcommand\Qcal{{\mathcal{Q}}}
\newcommand\Zcal{{\mathcal{Z}}}

\newcommand\eR{{\mathds{R}}}

\DeclareBoldMathCommand\oneb{1}
\DeclareBoldMathCommand\zerob{0}
\newtheorem{example}{Example}

\usepackage{cleveref}
\usepackage{float}
\usepackage{makecell}
\usepackage{arydshln}
\usepackage{booktabs}

\renewcommand{\phi}{\varphi}

\title{\LARGE \bf Adaptive Strategies for Pension Fund Management}

\author{Raphael Chinchilla$^{1*}$, Thomas D. Rueter$^{1}$, Timothy R. McDade$^{1}$, Peter R. Fisher$^{1}$,\\ Mykel J. Kochenderfer$^{1,2}$, Emmanuel Candès$^{1,2}$, Trevor Hastie$^{1,2}$, and Stephen Boyd$^{1,2}$% <-this % stops a space
  \thanks{\textsuperscript{*} Corresponding author: \href{mailto:raphael.chinchilla@blackrock.com}{raphael.chinchilla@blackrock.com} } % <-this % stops a space
  \thanks{$^{1}$BlackRock, 50 Hudson Yards, New York, NY 10001, USA.}%
  \thanks{$^{2}$Stanford University, Stanford, CA 94305, USA.}%
  \thanks{$\dagger$ The authors are employees of BlackRock and completed this as part of their employment there.}% <-this % stops a space
  \thanks{$\ddagger$ This paper represents the views of the authors and not of BlackRock, Inc.}% <-this % stops a space
}

\begin{document}
\maketitle
\thispagestyle{empty}
\pagestyle{empty}

\begin{abstract}
  This paper proposes a simulation-based framework for assessing and improving the performance of a pension fund management scheme. This framework is modular and allows the definition of customized performance metrics that are used to assess and iteratively improve asset and liability management policies. We illustrate our framework with a simple example that showcases the power of including adaptable features. We show that it is possible to dissipate longevity and volatility risks by permitting adaptability in asset allocation and payout levels. The numerical results show that by including a small amount of flexibility, there can be a substantial reduction in the cost to run the pension plan as well as a substantial decrease in the probability of defaulting.
\end{abstract}


\section{Introduction}

Retirement savers and investors confront a variety of risks, including market volatility, inflation, fund management cost, and longevity. Individuals whose retirement is funded by Defined Contribution (DC) schemes (such as a 401(k) in the United States) rather than Defined Benefit (DB) (as a pension) must bear these four risks themselves. For most retirement savers, the risk of outliving accrued savings is not negligible.

The framework we present in this paper is applicable to a variety of retirement financing arrangements. Although it was motivated by considering a DB plan in runoff, the framework today presents a more compelling case for the DC space for two reasons. First, DC has become a far more common arrangement than DB. Second, in DC retirement plans, no beneficiary has a guaranteed floor on payouts for the duration of their lifetime; any floor on lifetime payouts, like the one our framework offers, is a great value indeed. While this framework might seem standard in the control community, to the authors' best knowledge, this is the first time it is applied to a pension plan.

In DB pensions, beneficiaries are insensitive to capital market returns, inflation, fund management fees, or longevity risk. The corporate sponsor absorbs these risks, guaranteeing the beneficiary a lifetime cash flow regardless of the behavior of stocks, bonds, currencies, interest rates, inflation, fund management fees, or longevity patterns. In DC schemes, these risks are borne by the individual investor.

Longevity risk remains mostly unmanaged. Pension funds that confront mounting financial uncertainty, including that stemming from longevity risk, often consider undertaking a Pension Risk Transfer (PRT) where the pension plan transfers their assets and liabilities (payments to beneficiaries) to an insurance company, who then becomes responsible for paying their beneficiaries. The pension plan must pay a steep premium for the privilege of offloading this risk. Regulators require the insurers to hold an extremely conservative bond portfolio while holding large amounts of capital (upon which their shareholders expect appropriate returns). The rigidity of this structure is designed to protect the promises to the beneficiary but ends up being costly.

Longevity risk is salient to individual investors who are increasingly funding their own retirement, as well as to those institutional investors who still hold longevity risk on their balance sheet. Instead of relying on either a corporate sponsor or government balance sheet, or on transferring the uncertainties to an insurer, schemes could be designed to absorb both longevity and volatility risks by being adaptable in both the assets held and in the level of payouts over time. In this way, adaptable schemes could provide higher expected returns while still largely maintaining predictability.

The key contribution of this paper is a theoretical framework to assess the performance of an adaptable scheme through simulation. We start by specifying a model that describes the yearly evolution of the pension's assets given the invested portfolio returns, the net liabilities, and external funding. These three inputs depend on choices of portfolio, liabilities' management, and calls for external funding. Instead of fixed actions, which is the industry standard, our simulator requires policies that determine which actions to take each year given the available information. The policies are then evaluated by simulating them over randomly generated scenarios and collecting metrics about their performance. These metrics are quantities relevant to the stakeholders, such as default probability and investment return. If the performance metrics are not satisfactory, then the policies are updated and re-simulated, iteratively, until an acceptable policy is found. This framework is explained in \cref{sc:abstractmodel}.

\Cref{sc:ourimplementation} outlines a simple example of the abstract framework. The portfolios are constructed using a long-only Markowitz approach \cite{boyd2024markowitz} in which we include holding limits on certain asset classes. The liabilities are obtained from stochastic simulations of individual mortality, contributions paid into the fund, and pensions paid out. To manage these liabilities, the pensions paid out can be increased or decreased with respect to a baseline value. For the policy, we propose one that depends only on the value of the current assets under management divided by the total future (predicted) liabilities, commonly known as the asset to liability ratio. We discretize the range of this ratio into bins, and we assign which actions to take from each bin. Finally, we propose an automatic method to optimize the policy with respect to the realized metrics obtained from simulation.

In \cref{sc:numericalresults}, our numerical results show that by including a small amount of flexibility, there is a substantial reduction in the cost to run the pension plan as well as a substantial decrease in the probability of defaulting.

\subsection*{Related work}

\paragraph{Life cycle investing} Economics literature has long studied an individual’s ability to smooth consumption throughout their lifetime. Early work by Friedman et al. \cite{friedman1957permanent} on “permanent income” and Modigliani et al. \cite{modigliani1954utility} on the “life cycle” argued that individuals prefer stable consumption over their lifetime. This body of work has evolved into a robust analysis of multi-period optimization problems that individual economic actors can undertake, and the financial decisions that those processes imply \cite{cocco2005consumption}. Complex models in this field of study consider the way an individual balances current and future utility of consumption with mortality probabilities, a desire to bequeath assets to heirs, and risk aversion.

\paragraph{Annuities} Humans are not rational economic actors and do not perfectly smooth their consumption over the long run \cite{tirole2017economics}. Many families do not have stable income or consumption characteristics over time \cite{morduch2017financial}. One way to help ensure consumption remains stable despite shocks is purchasing insurance, particularly an income-oriented insurance product such as an annuity. Annuities pool mortality risks among beneficiaries, so those that live shorter than expected end up receiving an income stream equivalent to less than their initial contribution, leaving the remainder to subsidize those who outlive expectations. Annuity providers  receive considerable regulatory scrutiny and constrains the set of permissible assets an insurer can hold against their liabilities. As a result of the stability that annuities promise, their cost is often viewed as high. Scholars have studied annuity pricing over recent decades \cite{bauer2010pricing} and the social factors influencing demand for them \cite{boyer2020demand}.

\paragraph{Mortality pooling} Because of the variety in household economic circumstances, there is an extensive literature considering the best way for households to balance portfolio investment with annuity purchases \cite{pang2010optimizing}. One important subset of this literature draws inspiration from history to analyze asset management solutions to lifetime income using mortality pooling, some of which are called tontines \cite{goldsticker2007mutual, forman2014tontine, milevsky2015king, fullmer2025minimizing}. The basic idea of these solutions is that individuals contribute to a fund which is invested over time, and pay the members an income during the retirement phase. When a member of the fund dies, some the funds they contributed remain in the pool, increasing the payout of the remaining members. The mortality pooling concept has been proposed for the modern-day financial services industry \cite{milevsky2022build}. In recent years, the industry has evolved different mechanisms to distribute mortality credits in an actuarially fair manner, even across a pool of heterogeneous individuals \cite{fullmer2019individual}.

\section{Pension fund model} \label{sc:abstractmodel}

\subsection{Dynamics of a pension fund} \label{sc:abstractdynamics}

Let $v_t$ be the total value of the assets of a pension plan at the beginning of year $t$. Our goal is to create a model of the value of the assets of the pension plan at the beginning of the next year. This model depends on four components.

The first component is the return of the investments. We assume a universe of $n$ investible assets, denoting their returns during year $t$ by $\alpha_t \in \eR^n$.  Given a choice of portfolio $p_t$, the portfolio return is given by $r(\alpha_t,p_t)$. The portfolio return may depend non-linearly on the assets' returns $\alpha_t$.

The next component is the liabilities. Let $\iota_t$ be attributes (such as age and income) of an individual enrolled in the pension plan during year $t$. Given a plan's operational decisions $q_t$ for year $t$, the cashflow from an individual is denoted by $\lambda(\iota_t,q_t)$. The cashflow is positive if the plan is collecting premia from the individual and negative if the plan is paying the individual a pension. We denote by $\pi_t$ the collection of the attributes of all individuals in the plan during year $t$. Then, the total net liabilities of the plan is $\ell(\pi_t,q_t):=\sum_{\iota_t \in \pi_t} \lambda(\iota_t,q_t)$.

The last two components are more straightforward. Some external funding may be injected into the plan in year $t$, which we denote $e_t$. Finally, third parties may charge a percentage annual fee $m$, to manage the assets.

The total value of the assets of the plan at the beginning of the year $t+1$ is
\begin{equation}\label{eq:dynamics}
  v_{t+1}=(1+r(\alpha_t,p_t))(v_t + \ell_t(\pi_t,q_t) + e_t)(1-m).
\end{equation}
This model is agnostic on whether the values are expressed in nominal or real amounts (i.e. adjusted for inflation). An inflation adjustment would be reflected in the values of $r(\alpha_t,p_t)$, $\ell_t(\pi_t,q_t)$ and $e_t$ while taking the value of $v_0$ as reference. Moreover, this model assumes negligible transaction costs, but they could be included as an extra term.

\subsection{Pension policy/decision strategy} \label{sc:abstractpolicy}

At the beginning of each year, the pension plan's manager can take the following actions:
\begin{enumerate}
  \item Change the portfolio $p_t$ from a set of possible portfolios $\Pcal$.
  \item Ask the sponsor of the plan for an injection of external cash $e_t \in \eR_{\ge 0}$.
  \item Select an operational decision $q_t$ that modifies the net liabilities. The set of permissible options $\Qcal$ is determined by the pension's rules. These options may include modification to premia and/or payout, or offloading of liabilities via pension risk transfer.
\end{enumerate}
For convenience, we will denote all the actions taken at time $t$ by $a_t$, and the set of possible actions $\Acal \subset \Pcal\times\eR_{\ge 0}\times \Qcal$.

If there were no unknowns quantities, at $t=0$ the pension plan manager could specify which actions to take at time $t=1,2,\dots$ and guarantee that all of its liabilities would be met without ever risking default. In reality, pension plans manager react to the changing environment. Using a policy will allow us to model this recourse.

Let $\Zcal$ be the set of all possible information available at any time. We denote by $z_t\in \Zcal$ the specific information available at time $t$. In general, $z_t$ will contain information both about known and unknown quantities. Examples of known quantities are management fees $m$, current and past total value of assets $v_t,v_{t-1},\dots$, previous asset returns $\alpha_{t-1},\alpha_{t-2},\dots$ and previous population attributes $\pi_{t-1},\pi_{t-2},\dots$. The information for unknown quantities (present or future), may either be an estimate or a probability distribution that models their behavior. Examples of unknown quantities are future returns $\alpha_{t},\alpha_{t+1},\dots$ and future population attributes $\pi_{t},\pi_{t+1},\dots$ A policy $\phi_\theta(\cdot)$ is a function parameterized by $\theta$ such that $\phi_\theta(z_t)=a_t$. In other words, it determines which actions to take at each time given the available information with some choice of behavior captured in its parameters.

Before moving on, it is worth noting a few attributes of our framework. First, our approach involves changing actuarial liabilities and assets simultaneously. This is not the present approach of the industry, where current practice is to model actuarial liabilities first, then convert them into contractual liabilities, against which assets are managed. Second, it is possible to apply our simple framework in a more complicated mortality pooling context \cite{fullmer2019individual}. Third, as a contingent asset, a plan could also shed longevity or volatility risk for a price, for a period of time---as opposed to until the last beneficiary dies, at any price.

\subsection{Simulating and evaluating a policy}

To assess the effectiveness of a policy $\phi_\theta(\cdot)$ parameterized by $\theta$, we simulate the values of $v_t$ according to \eqref{eq:dynamics} for $t=1, 2, \dots, T$ and use these simulations to compute various performance metrics. We start by selecting a probability distribution to model the unknown variables $\alpha_{t}, \alpha_{t+1}, \dots$ and $\pi_{t}, \pi_{t+1}, \dots$. Using this distribution, we generate $S$ sample trajectories. These trajectories allow us to derive sample paths for $v_t$ under the policy $\phi_\theta(\cdot)$. With these sample paths, we can then calculate key metrics. For example, a pension plan manager might be interested in the expected value of the assets, the empirical probability of defaulting, and the average payout. We denote the metrics for a policy with parameters $\theta$ by $m_\theta \in \eR^M$. This process is summarized by \Cref{alg:evalpolicy}.

\begin{algorithm}[H]
  \caption{Evaluating a policy using Monte Carlo}
  \label{alg:evalpolicy}
  \begin{algorithmic}[1]
    \State Start all world scenarios with an asset $v_0$
    \State Sample $S$ trajectories of asset returns $\alpha_{1}, \dots \alpha_{T}$ and population $\pi_{1}, \dots, \pi_{T}$
    \For{$s=1,\dots,S$}
    \For{$t=1,\dots,T$}
    \State Determine all the information available $z_t^{(s)}$
    \State Evaluate $\phi_\theta(z_t^{(s)})$ to obtain the actions $a_t^{(s)}=(p_t^{(s)},e_t^{(s)},q_t^{(s)})$.
    \State Sample from the distributions to obtain $r(\alpha_t^{(s)},p_t^{(s)})$ and $\ell(\pi_t^{(s)},q_t^{(s)})$.
    \State Compute $v_{t+1}^{(s)}=\Big(1+r(\alpha_t^{(s)},p_t^{(s)})\Big)\Big(v_t^{(s)} + \ell(\pi_t^{(s)},q_t^{(s)}) + e_t^{(s)}\Big)(1-m)$.
    \EndFor
    \EndFor
    \State Compute metrics
  \end{algorithmic}
\end{algorithm}

\subsection{Tuning the policy}

Our objective in tuning the policy is to achieve acceptable outcomes according to the performance metrics. Since there are multiple performance metrics (or equivalently, we model the performance as a single vector), we introduce a cost function $h:\eR^m \to \eR\cup \{\pm \infty\}$ that generates a relative priority for improving different metrics. Tuning the policy means finding the set of parameters $\theta$ that finds a (local) minimum of $h(m_\theta)$. After the tuning is finished, the policy is evaluated in a new set of samples to establish its out-of-sample performance. Examples of algorithms to tune the policy can be found in \cite{kochenderfer2015decision,kochenderfer2022algorithms}.

\begin{figure}[t]
  \centering
  \includegraphics[width=1.\linewidth]{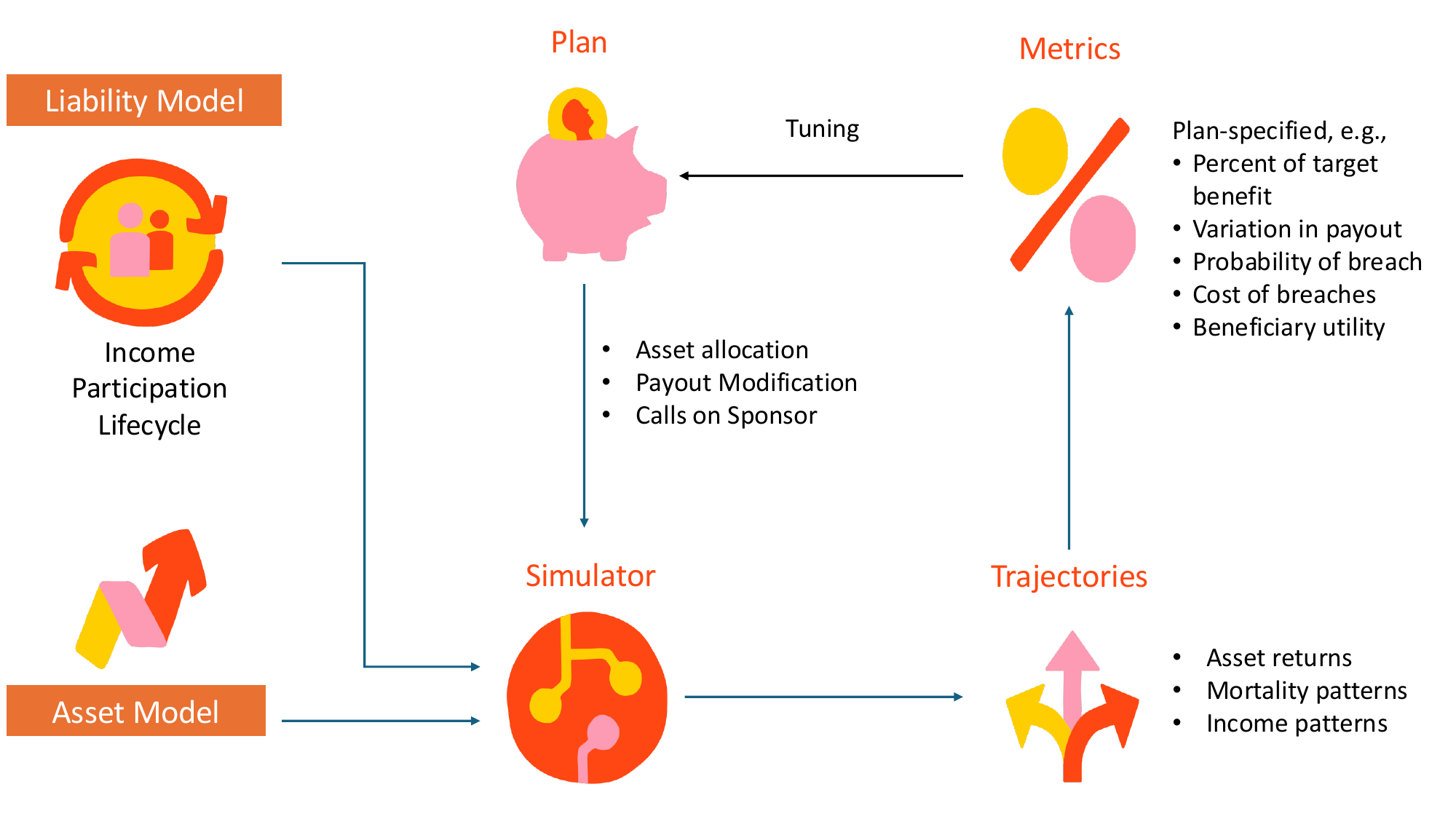}
  \caption{Schematic representation of the pension fund model.}
  \label{fig:schematic}
\end{figure}

\section{Example} \label{sc:ourimplementation}

The pension fund model of \cref{sc:abstractmodel} describes a general and abstract framework. In this section, we describe a specific example that we developed to illustrate these concepts. This is a notional example lacking some features that one would expect in an actual pension plan, such as using the duration and convexity as an input for the policy, but it is useful to show how even a simple example can lead to interesting results.

\subsection{Model of portfolio returns} \label{sc:modelofreturns}

Our model of returns is given by
\begin{equation}
  r(\alpha_t,p_t)=\alpha_t'w(p_t)
\end{equation}
where $\alpha_t$ are returns of asset classes and $w_t(p_t)$ are the weights associated to portfolio $p_t$.

Let $\bar \alpha$ be the expected return of the asset classes and $C$ the covariance matrix. We construct the weights of the portfolios $w(p_t)$ by using a Markowitz approach with risk aversion parameter $\mu_{p_{t}}>0$ \cite{boyd2024markowitz}. We define $w(p_t)$ to be the unique solution of the optimization problem
\begin{equation}\label{eq:markowitz}
  \begin{array}{ll}
    \mbox{maximize} &  \bar \alpha' w - \frac{\mu_{p_t} }{2} w'Cw \\
    \mbox{subject to}  &  \oneb'w=1\\
    & w \ge0 \\
    & w \le b_{in,up} \\
    & \sum_{i\in \text{domestic} }w_i \in [b_{do,low},b_{do,up}] \\
    & \sum_{i\in \text{foreign} }w_i \in [b_{fo,low},b_{fo,up}] \\
    & \sum_{i\in \text{fixed income} }w_i \in [b_{fi,low},b_{fi,up}] \\
    & \sum_{i\in \text{equities} }w_i \in [b_{eq,low},b_{eq,up}] \\
    & \sum_{i\in \text{alternatives} }w_i \in [b_{alt,low},b_{alt,up}]
  \end{array}
\end{equation}
where $b_{in,up}$ is the maximum allocation for any asset class and $[b_{do,low},b_{do,up}]$, $[b_{fo,low},b_{fo,up}]$, $[b_{fi,low},b_{fi,up}]$, $[b_{eq,low},b_{eq,up}]$ and $[b_{alt,low},b_{alt,up}]$ are lower and upper bounds on the total weights for asset classes that are domestic, foreign, fixed income, equities and alternatives, respectively.

\subsection{Model of liabilities} \label{sc:modelofliabilities}

As explained in \cref{sc:abstractdynamics}, the total liability is a sum over the cashflows of all the individuals in the pension plan. In our model, we simulate the life of each individual, or, more precisely, the change in the attributes $\iota_t$ through the life of each individual and when the individual dies. The individual's attributes we consider are biological factors (year of birth, biological sex) and socio-economic factors (income, place of residence, level of education). These attributes are used to construct an individualized mortality rate; we refer the reader to Appendix \ref{ap:coxmodel} for an explanation. Given each individual's mortality rate, we build a probability distribution on the trajectories of this individual's life over the next $T$ years. The individual cashflows $\lambda(\iota_t,q_t)$ is a function of these trajectories.

\paragraph{Sampling life trajectories}
During the individual's life over the next $T$ years, they can be in one of four states: alive and not yet working, alive and working, alive and retired, and dead. The probability of an individual dying is given by their mortality rate that year. If the individual is not dead, they can change state. If the individual is not yet working and is older than a minimum working age $m_w$, each year they have a probability of $1/y_w$ of starting to work ($y_w$ is the average number of years it would take for them to start working). If the individual is working and is older than a minimum retiring age $m_r$, each year they have a probability of $1/y_r$ of retiring ($y_r$ is the average number of years it would take for them to retire). Our model does not take into account unemployment, but this could easily be remediated by creating an extra state.

\paragraph{Individual cashflows} During an individual life trajectory, if they are working, they pay $k_w\%$ of their income; if they are retired, they receive $k_r\%$ of the last salary they received.

\paragraph{Plan rules and operational strategies} The set of permissible liability options $\Qcal$ (defined in \cref{sc:abstractpolicy}) is composed of $m_w$, $m_r$, $y_w$, $y_r$, $k_w$, and $k_r$. The values of $m_w$, $m_r$, $y_w$, $y_r$, and $k_w$ are determined by the plan rules and cannot be changed from one year to the other. The value of $k_r$ can be changed each year up to $\pm \tilde k_r\%$ from a baseline which we denote by $\bar k_r$ so $k_r=(100\%\pm \tilde k_r\%)\bar k_r$.

\begin{example}[Example of rules for a pension plan]
  The minimum age for starting to work is $m_w =18$ years old, with an average of $y_w = 4$ years to start doing so. When working, someone contributes $k_w\%=15\%$ of their income. The minimum age for retiring is $m_r = 60$ years old, with an average of $y_r=2$ years for doing so. When retired, someone is promised a baseline $\bar k_r = 80\%$ of their salary before retiring. The actual value of the pension each year will vary by up to $\pm \tilde k_r\% = \pm 10\%$ meaning that a pensioner is guaranteed to receive $(100\%\pm 10\%)80\%$ of their salary before retiring.
\end{example}

We reference the net liabilities according to a difference in the payout level $(100\%\pm \tilde k_r\%)$. For instance, $\ell(\pi_t,100\%)$ denotes the net liability when the payments to retirees are at the baseline level.

\subsection{Policy}

The information available at time $z_t$ is composed of:
\begin{enumerate}
  \item the current portfolio value $v_t$,
  \item the current baseline liabilities $\ell(\pi_t,100\%)$,
  \item the expected return of the asset classes defined by $\bar \alpha_t:=\mathbb{E}[\alpha_t]$,
  \item the expected net baseline liability for the next $T$ years given the population of year $t$, defined by $\bar \ell_{t+\tau}:=\mathbb{E}[\ell(\pi_{t+\tau},100\%)],\ \tau = 1,\dots,T$ provided by the model described in \cref{sc:modelofliabilities},
  \item the expected yield curve over the next $T$ years as of year $t$, denoted by $\gamma_{t+\tau},\ \tau = 1,\dots,T$.
\end{enumerate}

Let $L_t=\sum_{\tau=0}^{T}\gamma_{t+\tau} \bar \ell_{t+\tau}$ be the expected total net liabilities over the next $T$ years. The first step of the policy is to compute an extended version of the asset to liability ratio, defined by
\begin{equation*}
  \rho_t=
  \begin{cases}
    v_t/L_t & \text{if } v_t\ge0 \text{ and } L_t>0\\
    +\infty & \text{if } v_t\ge0 \text{ and } L_t\le0\\
    0 & \text{if } v_t<0
  \end{cases}\quad.
\end{equation*}
This formulation is motivated by the following observation. On one hand, if the total assets are negative, then the plan is in a dire state, independent of the future liabilities. On the other hand, when the assets are positive, if the future liabilities are not positive, then the plan is in a very good state.

We have already described in \cref{sc:modelofreturns,sc:modelofliabilities} the choices of portfolio $p_t$ and liability options $q_t$. The last action we need to describe is how $e_t$ is determined. The policy specifies a target asset to liability ratio at time $t+1$ denoted by $\widehat{\rho}_{t+1}$ and the external cash $e_t$ is given by
\begin{equation}
  \max\bigg\{\frac{\widehat{\rho}_{t+1}L_{t+1}}{(1-m)(1+w(p_t)'\bar \alpha_t)}-(v_t+\ell_t(q_t)),0\bigg\}
\end{equation}
which is the value of $e_t$ that would achieve an asset to liability ratio equal to $\widehat{\rho}_{t+1}$ at the beginning of year $t+1$ if the future liabilities and returns match the expected ones.

The policy is specified by discretizing the values of $\rho_t$ and defining which actions to take in each bin. This leads to a table in which each row specifies the actions to be taken at time $t$. \Cref{table:examplepolicy} is an example of what such table would look like. The parameters $\theta$ that indexes the policy are the values of the actions for each asset to liability ratio bin. The set of possible parameters $\Theta$ are the possible values of portfolio volatility, payout level and target asset to liability ratio.

\begin{table}
  \centering
  \begin{tabular}{@{}cccc@{}}
    \toprule
    \makecell{\textbf{Asset to} \\ \textbf{liability ratio}} & \makecell{\textbf{Portfolio} $p_t$ \\ \textbf{vol. (\%)}} & \makecell{\textbf{Cashflow option}  \\ $q_t$ \textbf{payout (\%)}} & \makecell{\textbf{Ext.} \\ \textbf{cash}} \\
    \midrule
    $> 2.0$ & 4.437 & 110 & 0 \\
    $1.8$--$2.0$ & 7.891 & 108 & 0 \\
    $1.6$--$1.8$ & 7.891 & 107 & 0 \\
    $1.4$--$1.6$ & 7.891 & 105 & 0 \\
    $1.2$--$1.4$ & 5.620 & 102 & 0 \\
    $1.0$--$1.2$ & 5.620 & 100 & 0 \\
    $0.8$--$1.0$ & 4.437 & 97 & 1.2 \\
    $< 0.8$ & 4.37 & 95 & 1.0 \\
    \bottomrule
  \end{tabular}
  \caption{Policy table: actions by asset-to-liability ratio}
  \label{table:examplepolicy}
\end{table}

\subsection{Metrics and policy tuning} \label{sc:PolicyTuning}
\paragraph{Policy evaluation} The policy is evaluated according to the following metrics:
\begin{itemize}
  \item yearly probability of performing a cash call $c$
  \item average payout level $\bar q$
  \item average change in payout level from one year to the next $\Delta q$
\end{itemize}

We scalarize the metrics using the following function. For each metric value $m$ we denote its lower acceptable bound by $m^{\text{low}}$, its upper acceptable bound by $m^{\text{up}}$, its midpoint by $m^{\text{mid}} = (m^{\text{low}}+m^{\text{up}})/2$, and its priority (weight) by $m^{\text{pri}}$. We also define the inner slope scale $s_m = 2/(m^{\text{up}}-m^{\text{low}})$ and use $10 s_m$ as the outer slope (penalty amplification outside the acceptable band). Then for each metric $m$ we associate the value
\begin{align*}
  \eta(m)=
  \begin{cases}
    10 m^{\text{pri}} s_m ( m^{\text{low}} - m )+1 & m < m^{\text{low}} \\
    m^{\text{pri}} s_m ( m^{\text{low}} - m )+1 & m^{\text{low}} \le m < m^{\text{mid}} \\
    m^{\text{pri}} s_m ( m - m^{\text{up}} )+1 & m^{\text{mid}} \le m < m^{\text{up}} \\
    10 m^{\text{pri}} s_m ( m - m^{\text{up}} )+1 & m \ge m^{\text{up}}
  \end{cases}
\end{align*}
whose plot can be seen in \Cref{fig:penalty_function}.
\begin{figure}[t]
  \centering
  \begin{tikzpicture}
    \begin{axis}[
        width=8cm,
        height=4cm,
        xlabel={$m$},
        ylabel={$\eta(m)$},
        grid=major,
        grid style={gray!30},
        axis lines=center,
        xmin=-0.5,
        xmax=3.5,
        ymin=0,
        ymax=9,
        xtick={1, 1.5, 2},
        xticklabels={$m^{\mathrm{low}}$, $m^{\mathrm{mid}}$, $m^{\mathrm{up}}$},
        ytick={0},
        thick,
      ]

      % Define parameters for the plot
      \def\mlow{1}
      \def\mup{2}
      \def\mmid{1.5}
      \def\mpri{1}
      \def\sm{0.5}

      % Region 1: m < m^low (steep negative slope)
      \addplot[blue, very thick, domain=0:1] {20 * \mpri * \sm * (\mlow - x) + 0.25};

      % Region 2: m^low ≤ m < m^mid (moderate negative slope)
      \addplot[blue, very thick, domain=1:1.5] {\mpri * \sm * (\mlow - x) + 0.25};

      % Region 3: m^mid ≤ m < m^up (moderate positive slope)
      \addplot[blue, very thick, domain=1.5:2] {\mpri * \sm * (x - \mup) + 0.25};

      % Region 4: m ≥ m^up (steep positive slope)
      \addplot[blue, very thick, domain=2:3] {20 * \mpri * \sm * (x - \mup) + 0.25};

      % Add vertical lines at key points
      \addplot[dashed, gray] coordinates {(1, 0) (1, 9)};
      \addplot[dashed, gray] coordinates {(1.5, 0) (1.5, 9)};
      \addplot[dashed, gray] coordinates {(2, 0) (2, 9)};

      % Add points at transitions
      \addplot[mark=*, mark size=3pt, only marks, black] coordinates {(1, 0.25) (1.5, 0) (2, 0.25)};

    \end{axis}
  \end{tikzpicture}
  \caption{The penalty function $\eta(m)$ with steep penalties outside the acceptable and moderate penalties within it.}
  \label{fig:penalty_function}
\end{figure}

The objective function that is tuned is $h(m_\theta)=\eta(c)+\eta(\bar q)+\eta(\Delta q)$.

\paragraph{Direct policy optimization} We would like to find the optimal parameters $\theta$ that achieve the best cost $h(m_\theta)$. We do this using a derivative-free algorithm that achieves one-optimality, meaning that the policy cannot be improved by only changing the value of the action in one cell. We start by creating a discrete set of possible values of portfolio volatility, of payout levels, and of target asset to lability ratio $\widehat{\rho}$.
\begin{algorithm}[H]
  \caption{Algorithm to achieve one-optimality}
  \label{alg:oneopt}
  \begin{algorithmic}[1]
    \State Initialize table with values that seem reasonable.
    \State Start all world scenarios with an asset $a_0$
    \State $\text{counter}=0$
    \While{counter $\le$ number of cells in the table}
    \For{each row of the table}
    \For{each column of the table}
    \For{each allowed action corresponding to the current column}
    \State Evaluate the policy using this action for this cell
    \EndFor
    \If{any of the tested action produces a smaller value than the current best action}
    \State Assign this new best action for this cell
    \State $\text{counter}=0$
    \Else
    \State $\text{counter}=\text{counter}+1$
    \EndIf
    \EndFor
    \EndFor
    \EndWhile
  \end{algorithmic}
\end{algorithm}

\section{Results} \label{sc:numericalresults}

To demonstrate our framework, we choose implementations of the framework that take place in the United States. For the purpose of this research study, we use the demographic data from the Health and Retirement Study from the University of Michigan \cite{rand2024longitudinal}. The HRS (Health and Retirement Study) is sponsored by the National Institute on Aging (grant numbers NIA U01AG009740 and NIA R01AG073289) and is conducted by the University of Michigan \cite{hrs2020}. This is a nationally representative data set that includes extensive information about individuals sociodemographic and economic characteristics, including income and mortality. From this dataset, we use the information from income, level of education, residency, biological sex, data of birth and date of death which we combine with the mortality rate projections created by the United Nations Department of Economic and Social Affairs Population Division \cite{un_population_prospects_2024} in order to create a Cox proportional hazards model, as described in Appendix \ref{ap:coxmodel}.

% RC: Emphasizing again that the dataset from the rand corp is for research only

The asset classes returns are generated using a proprietary tool based on the work of Van Beek \cite{vanbeek2020consistentcalibrationeconomicscenario}. It generates random trajectories of annual returns for asset classes and macroeconomic projections for over 600 indicators for 50 years. Of these indicators, we selected a set of 25 asset classes denominated in U.S. Dollars. They include bonds (\eg U.S. Treasury bonds with maturities of 1, 3, 5, 7, 10, 20+ years), equities (\eg U.S. large cap, U.S. small cap, China A Shares) and some alternatives (\eg infrastructure debt, real estate).

The characteristics of the implementation are:
\begin{itemize}
  \item \textbf{Population:} The simulated population is composed of 1000 individuals with gender, ages and socio-economic background that are representative of the United States population. The plan is closed to new entrants that are born after the initial year, but is open to entrants that are already born.

  \item \textbf{Plan rules:} Each individual is allowed to start working when they are 18 years old, and, on average, start working when they are 22 years old. While they are working, they contribute 10\% of their income. They are allowed to retire once they turn 60 years old, and, on average, retire when they are 65 years old. Once retired, they expect a baseline pension equivalent of 80\% of their last income.

  \item \textbf{Simulation length:} The policy is simulated for 30 years, from 2025 to 2055. Each year, the total liabilities are computed by summing the projected liabilities over the following 30 years discounted by the projection of the yield curve at that date. So, for instance, in the year 2035, the total liabilities consist the projected liabilities from 2035 to 2065, discounted by the projected yield curve of 2035.

  \item \textbf{Initial capitalization:} All plans start with the same amount of assets. This initial amount of asset was chosen such that the plan start with an asset to liability ratio of $1.6$.

  \item \textbf{Asset to liability bins:} We discretize the asset to liability ratio into the following bins: $[0.0, 0.2)$, $[0.2, 0.4)$, $\dots$, $[2.2,\infty)$.
\end{itemize}

By running our machine with different pension plan design features, we can test the effectiveness of various pension plan set-ups.
Our ultimate aim is to be able to compare the efficiency of different pension plans across a spectrum from high cost (like a PRT-insurer provided annuity) to low cost (like a flexible “adaptable annuity” that dynamically changes asset allocation and payout level, subject to plan design constraints).
We run four different instantiations of pension plan designs within our framework.
\begin{itemize}
  \item \textbf{Plan A:} Designed to be an approximation of a Pension Risk Transfer annuity (fixed income only, static allocation, paying beneficiaries 100\% of target always). Various costs associated with a realistic PRT, such as the costs of hedging, cost of capital, or pension insurance premia, are not included, but would only further inflate the costs to operate this plan.
  \item \textbf{Plan B:} Introduces flexibility on the payouts: fixed income only, static allocation, paying beneficiaries between 90--110\% of target (with payout level varying by no more than 2\% of target in any year).
  \item \textbf{Plan C:} Introduces flexibility on the asset allocation as well: fixed income only, variable allocation, paying beneficiaries between 90--110\% of target (with payout level varying by no more than 2\% of target in any year).
  \item \textbf{Plan D:} Introduces equities in the asset allocation, in addition to fixed income.
\end{itemize}
\Cref{tab:plans} summarize the difference between the plans.

\begin{table}
  \centering
  \begin{tabular}{|l|c|c|c|c|}
    \hline
    Plan & \textbf{A} & \textbf{B} & \textbf{C} & \textbf{D} \\ \hline
    \textbf{Benefit} & & & & \\
    \% of target & 100\% & \makecell{90\%--\\110\%} & \makecell{90\%--\\110\%} & \makecell{90\%--\\110\%} \\
    \hdashline
    \makecell[l]{\textbf{Assets} \\ \textbf{available}} & \makecell{Fixed\\income} & \makecell{Fixed\\income} & \makecell{Fixed\\income} & \makecell{Fixed \\ income\\ + equities} \\
    \hdashline
    \makecell[l]{\textbf{Asset} \\ \textbf{allocation}} & \makecell{Static,\\rebalance} & \makecell{Static,\\rebalance} & Adaptive & Adaptive \\
    \hline
  \end{tabular}
  \caption{Differences between Plans A, B, C, and D}
  \label{tab:plans}
\end{table}

\begin{figure*}[tbp]
  \centering
  \begin{subfigure}[t]{0.24\textwidth}
    \centering
    \includegraphics[width=\linewidth,height=0.8\textwidth]{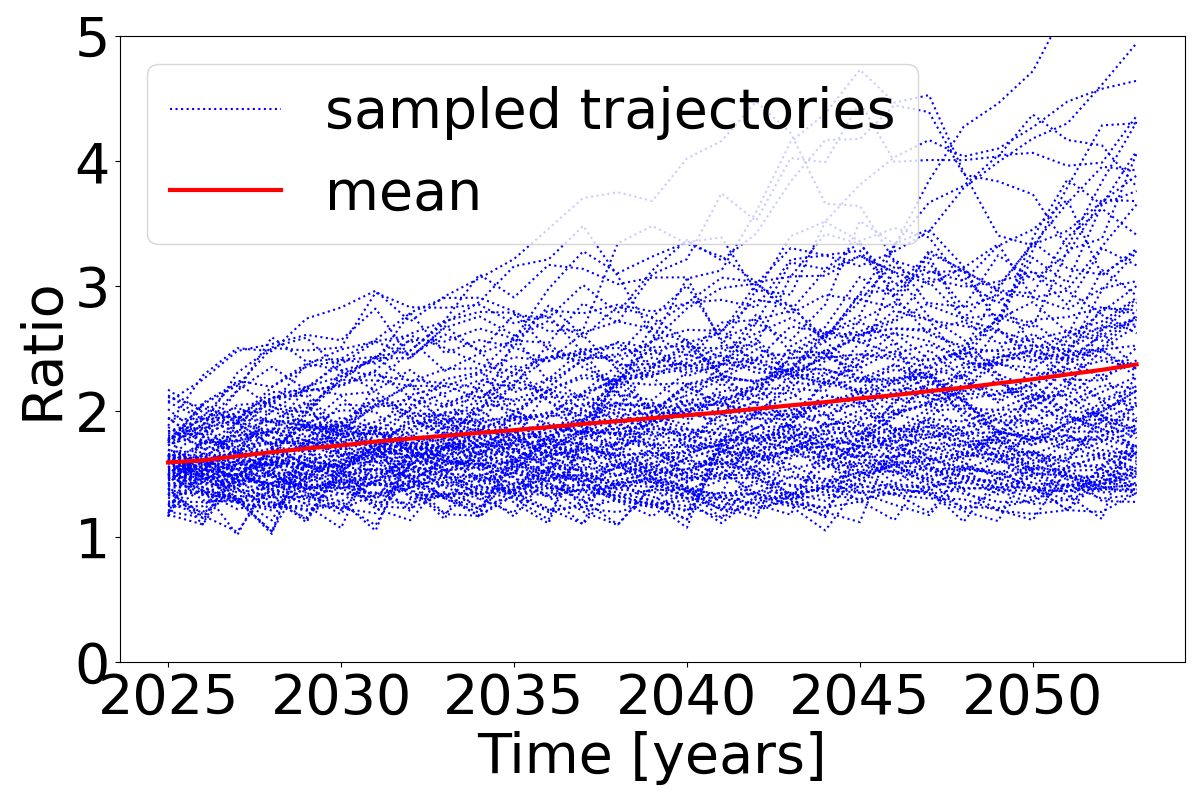}
    \caption{Plan A}
    \label{fig:results_A}
  \end{subfigure}
  \hfill
  \begin{subfigure}[t]{0.24\textwidth}
    \centering
    \includegraphics[width=\linewidth,height=0.8\textwidth]{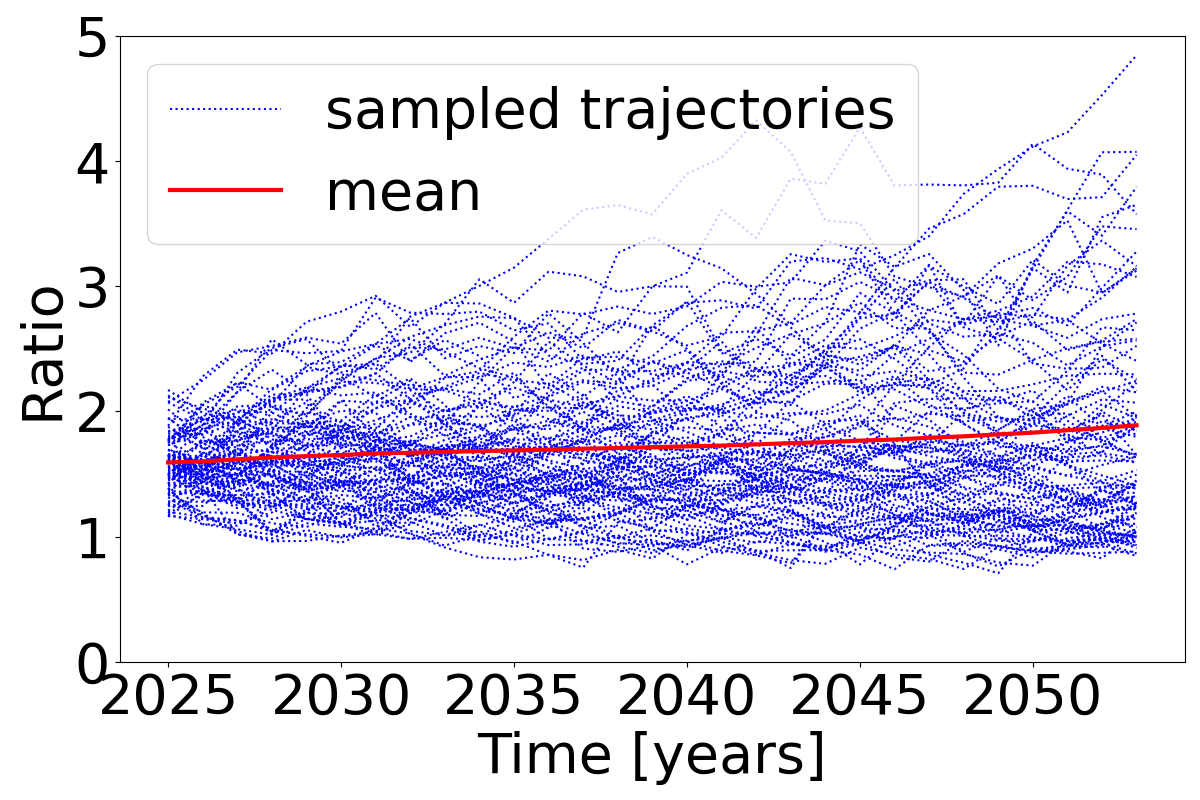}
    \caption{Plan B}
    \label{fig:results_B}
  \end{subfigure}
  \hfill
  \begin{subfigure}[t]{0.24\textwidth}
    \centering
    \includegraphics[width=\linewidth,height=0.8\textwidth]{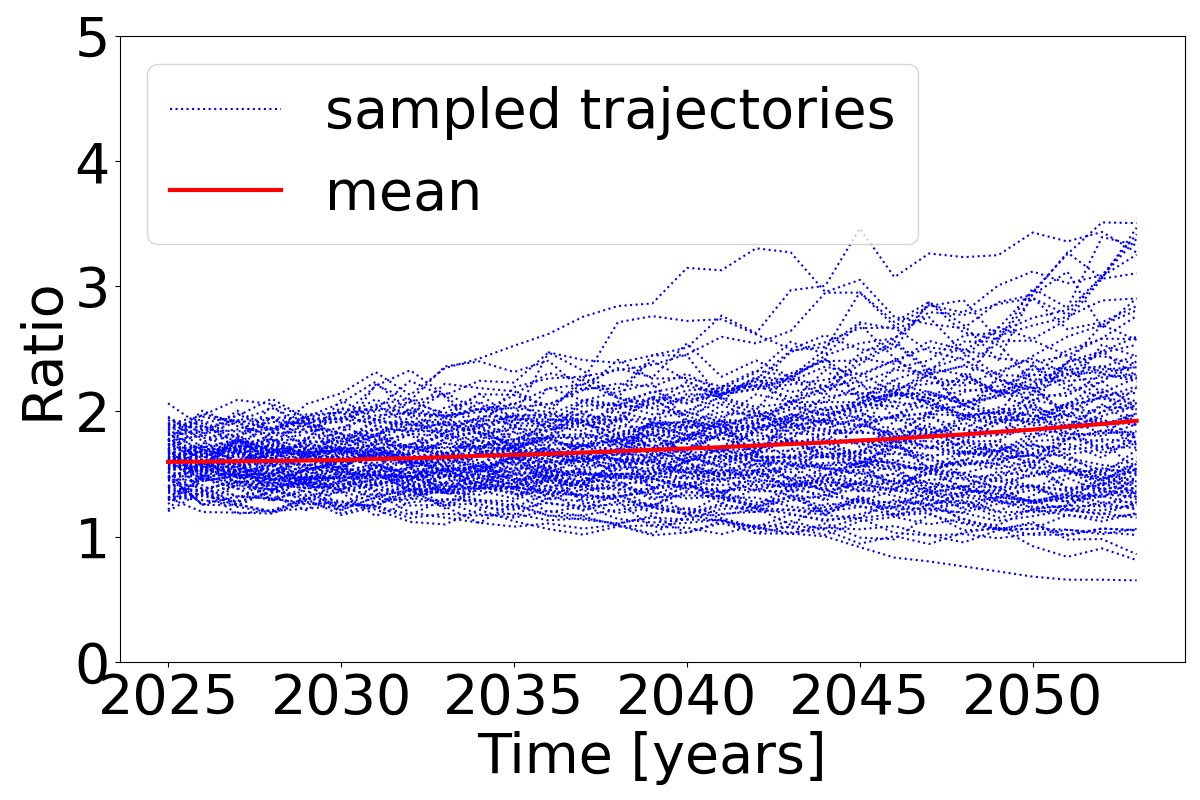}
    \caption{Plan C}
    \label{fig:results_C}
  \end{subfigure}
  \hfill
  \begin{subfigure}[t]{0.24\textwidth}
    \centering
    \includegraphics[width=\linewidth,height=0.8\textwidth]{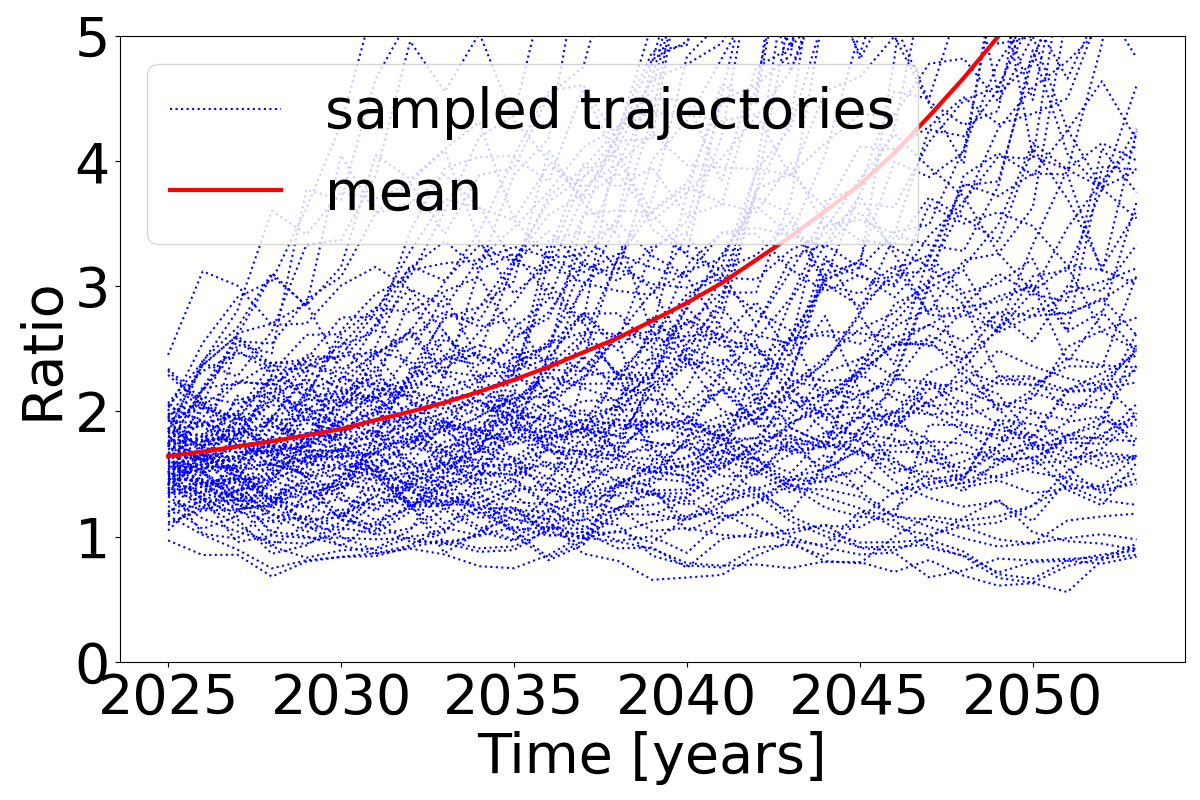}
    \caption{Plan D}
    \label{fig:results_D}
  \end{subfigure}
  \caption{Simulations of asset to liability ratio under the different plans.}
  \label{fig:results_all}
\end{figure*}

We run each plan over 10000 simulated world scenarios (\ie 10000 trajectories of liabilities and 10000 trajectories of asset returns) and compute metrics of these results. The results are available in \Cref{tab:comparison} and \Cref{fig:results_A,fig:results_B,fig:results_C,fig:results_D}. Increasing the degrees of freedom of the plans, improves performance, either by distributing a larger benefit to the pensioners or by decreasing the risk of requiring a cash call. Notice that the level at which each plan is considered breached is automatically learned, and is the optimal level given the plan's degrees of freedom. This happens even though all the plans share the same objectives. The automatically learned breach asset-to-liability thresholds are 1.2, 1.0, 0.6, and 0.6 for Plans A, B, C, and D, respectively.

\begin{table}
  \centering
  \begin{tabular}{@{}lcccc@{}}
    \toprule
    \textbf{Metric} & \textbf{A} & \textbf{B} & \textbf{C} & \textbf{D} \\ \midrule
    \makecell[l]{\textbf{Mean \% of} \\ \textbf{target benefit}} & 100\% & 99.9\% & 102.4\% & 104.0\% \\ \midrule
    \makecell[l]{\textbf{Prob. breaching} \\ \textbf{A:L floor 1y (30y)}} & \makecell{1.32\% \\ (53\%)} & \makecell{0.88\% \\ (25\%)} & \makecell{0.21\% \\ (7.7\%)} & \makecell{0.19\% \\ (6.0\%)} \\ \midrule
    \makecell[l]{\textbf{Ex post breach} \\ \textbf{value (\% initial)}} & 38\% & 13\% & 2.9\% & 3.0\%  \\
    \bottomrule
  \end{tabular}
  \caption{Performance metrics for Plans A--D}
  \label{tab:comparison}
\end{table}

It is worth reiterating that the Plan A described here, even with its relative expense, is still cheaper than an annuity would be in practice. Our Plan A does not deduct from returns the cost of capital or hedging risks from interest rates or foreign exchange fluctuations. When these costs are included, we expect that Plans B, C, and D will be even cheaper by comparison.

\section{Conclusion}

The principal contribution of this paper is a general framework to address both longevity risk and return uncertainty in providing a low(er) cost source of lifetime income. The core idea is to create a simulator that generates trajectories of asset returns and liabilities and use them to tune a policy. While this idea is widespread in the field of control, applying an engineering solution of this type appears to be less common in finance.

Every country, every pension system, every individual will ultimately provide “adaptive” answers to the problem of sustaining consumption in the future from whatever pot of assets may have been accumulated. By directly approaching the idea of an adaptive pension seems an important place to start.

Our second contribution was a specific policy. While this policy is rather simple, it was used to illustrate a key result: including just a small amount of flexibility can substantially decrease risk and cost by allowing the pension plan to adapt to new circumstances. In our simple policy, this was restricted to only changes in portfolio allocation and pension payout, but could be used in an array of other decisions that are currently fixed by plan rules. If designed correctly, this approach would benefit both pensioners and sponsors because it trades small amounts of risk for substantial decreases in cost and increase in reliability.

\appendix

\subsection{Mortality rates through the Cox model}\label{ap:coxmodel}

Each individual is characterized by three core characteristics: year of birth, biological sex, and country of residency. In addition to these core characteristics, individuals are also described by socio-economic covariate, such as income, level of education, etc. We assume we have at our disposal a dataset which contains individuals' covariates over a period of some years and their dates of deaths or whether they are alive by the end of the period.

There are several institutions around the world that provide a mortality rate for an individual given the three core characteristics. Our goal is to understand how this mortality rate is affected by covariates.

One way to model longevity for an individual with covariate $X \in \mathbb{R}^d$ is to use the Cox model \cite{cox1972regression}. The hazard function is given by
\begin{equation*}
  P(T \in [t, t + dt] \mid T \ge t,  X) = \lambda(t) \exp(\theta^\top X);
\end{equation*}
here,
\begin{enumerate}
  \item $\lambda(t)$ is the baseline mortality
  \item $\exp(\theta^\top X)$ is an exponential term that captures the effect of the covariates on the hazard function;
  \item $\theta \in \mathbb{R}^d$ is the vector of regression coefficients which needs to be estimated from the data.
\end{enumerate}

\paragraph{The survival function}
Let $T$ be the random variable representing the time of death. Let $f_\theta(T, X)$ be the density function of $T$ and $\bar{F}_\theta(T,X)$ be the survival function (the chance of making to at least $T$). Then by definition,
\begin{equation*}
  \frac{\text{d} \log \bar{F}_\theta(t, X)}{\text{d}t} = - \lambda(t) \exp (\theta^\top X).
\end{equation*}
This implies that the survival function is given by
\begin{equation*}
  \bar{F}_\theta(t, X) = \exp\left(-\int_{t_0}^t \lambda(s) \exp (\theta^\top X ) \, \text{d}s\right),
\end{equation*}
where $t_0$ represents the moment when the individual starts to be tracked. It also follows that the density function is given by
\begin{equation*}
  f_\theta(t, X) = \lambda(t) \exp(\theta^\top X) \, \bar{F}_\theta(t,X).
\end{equation*}

\paragraph{Maximum likelihood estimation}
We have i.i.d.~data $(T_i, X_i), \ldots, (T_n, X_n)$, where $T_i$ is either the time of death or the time of censoring for the $i$th individual and $X_i$ is the covariate vector for the $i$th individual. Let $\delta_i$ be the indicator of whether the time of death of the $i$th individual is observed. The likelihood function for the Cox model is given by
\begin{equation*}
  \prod_{i=1}^n \left( f_\theta(T_i,X_i) \right)^{\delta_i} \left( \bar{F}_\theta(T_i,X_i) \right)^{1-\delta_i}.
\end{equation*}
Therefore, the log-likelihood is given by $\ell(\theta) = $
\begin{equation*}
  \sum_{i=1}^n \left( \delta_i \left( \log \lambda(T_i) + \theta^\top X_i \right) - \exp(\theta^\top X_i)   \int_{t_0}^{T_i} \!\! \lambda(s) \, \text{d}s \right).
\end{equation*}
Setting $\Lambda(t) = \int_{t_0}^{t} \lambda(s) \, \text{d}s$, the maximum likelihood estimator $\hat{\theta}$ is the solution of
\begin{equation*}
  \min_\theta\ \sum_{i=1}^n    \Lambda(T_i) \, \exp (\theta^\top X_i ) -\delta_i \theta^\top X_i.
\end{equation*}
This is a convex optimization problem and can be solved using standard optimization algorithms. (Recall that $\lambda(t)$ and, therefore, $\Lambda(t)$ are known.)

In our implementation, we are interested in understanding how the baseline mortality rate for a population in the United States changes according to income, years of education and location. They are treated as categorical:
\begin{itemize}
  \item Income is separated in 4 quantiles (rich, upper middle class, lower middle class, poor) each representing 25\% of the population.
  \item Education is separated into 4 levels (high school dropout, high school degree, college degree, advanced degree).
  \item Location is separated into the 4 census regions (Northeast, South, Midwest, West).
\end{itemize}

\subsection{Model of income}

The income of an individual is determined by their age and the income bin in which they categorize (see the end of the previous section). Let $\alpha$ be the individual age. Then their income is given by
\begin{equation*}
  \text{income}(\alpha)=
  \begin{cases}
    25\beta & \text{if } \alpha< 25 \\
    \alpha\beta & \text{if } \alpha \in [25,45] \\
    45\beta & \text{if } \alpha> 45 \\
  \end{cases}
\end{equation*}
where $\beta$ is a constant that depends on the individual's income bin (which is fixed throughout a person's life). These values of $\beta$ are learned from our dataset. For our dataset, we obtained poor $=325$, lower middle class $=750$, upper middle class $=1400$, and rich $=2880$ (rounded to the three significant digits). The shape of this income function is illustrated below:
\begin{center}
  \begin{tikzpicture}[scale=0.6]
    \begin{axis}[
        width=8cm,
        height=5cm,
        xlabel={Age $\alpha$},
        ylabel={Income},
        xmin=20,
        xmax=60,
        ymin=0,
        ymax=50,
        xtick={25, 45, 60},
        ytick={25, 45},
        yticklabels={$25\beta$, $45\beta$},
        yticklabel style={font=\large},
        axis lines=left,
        grid=none,
      ]
      \addplot[blue, very thick, domain=20:25] {25};
      \addplot[blue, very thick, domain=25:45] {x};
      \addplot[blue, very thick, domain=45:60] {45};
      \addplot[blue, mark=*, mark size=2pt, only marks] coordinates {(25, 25) (45, 45)};
    \end{axis}
  \end{tikzpicture}
\end{center}

\bibliographystyle{ieeetr}
\bibliography{biblio}

%%%%%%%%%%%%%%%%%
\end{document}